\documentclass[12pt,leqno]{amsart}
\usepackage{amsmath,amsfonts,amsthm,amssymb,amscd, dsfont,tipa,mathtools,listings, rotating, multirow, bbm}
\usepackage{color}
\usepackage{comment}
\usepackage{seqsplit}
\usepackage{enumitem}
\usepackage{soul}
\usepackage{mathrsfs}
\usepackage[titletoc,title]{appendix}
\usepackage{tikz-cd}
\usepackage{lineno}
\usepackage{lipsum}
\usepackage{hyperref}
\usepackage{constants}
\newconstantfamily{abcon}{symbol=c}

\newtheorem{thm}{Theorem}[section]

\newtheorem{cor}[thm]{Corollary}

\newtheorem{rek}[thm]{Remark}

\DeclareMathOperator{\GL}{GL}

\DeclareMathOperator{\N}{\mathbb{N}}

\makeatletter
\let\@wraptoccontribs\wraptoccontribs
\makeatother

\title{The $\theta = \infty$ Conjecture and the Riemann Hypothesis for Automorphic $L$-functions}
\author{Anji Dong}
\address{Anji Dong: 
Department of Mathematics, University of Illinois Urbana-Champaign, Altgeld hall, 1409 W. Green Street, Urbana, IL, 61801, USA}
\email{anjid2@illinois.edu}

\author{Nawapan Wattanawanichkul}
\address{Nawapan Wattanawanichkul: Department of Mathematics, University of Illinois Urbana-Champaign, Altgeld hall, 1409 W. Green Street, Urbana, IL, 61801, USA}
\email{nawapanwattanawanichkul4@gmail.com}
    
\author{Alexandru Zaharescu}
\address{Alexandru Zaharescu: Department of Mathematics, University of Illinois Urbana-Champaign, Altgeld hall, 1409 W. Green Street, Urbana, IL, 61801, USA \and Simion Stoilow Institute of Mathematics  of the Romanian Academy, P.O. Box 1-764, RO-014700 Bucharest, Romania}
\email{zaharesc@illinois.edu}


\begin{document}

\begin{abstract}
The $\theta=\infty$ conjecture asserts that the mollified second moments of the Riemann zeta function remain bounded for mollifiers of arbitrary polynomial length. We investigate an analogue of this conjecture for automorphic $L$-functions associated with cuspidal representations of $\text{GL}_m(\mathbb{A}_{\mathbb{Q}})$, exploring its implications for the distribution of their nontrivial zeros. Extending the framework of Bettin and Gonek, we prove that if the mollified second moments of these $L$-functions remain suitably bounded for mollifiers of arbitrary polynomial length, then the $L$-functions are non-vanishing in corresponding regions of the critical strip. Furthermore, we establish a version of this criterion for families of $L$-functions, demonstrating that the $\theta = \infty$ conjecture for a family of $L$-functions implies a quasi-Riemann Hypothesis for that family.
\end{abstract}

\maketitle

\section{Introduction}
For $m\in\N$, let $\mathfrak{F}_m$ be the set of all cuspidal automorphic representations of $\text{GL}_m(\mathbb{A}_Q)$ whose central characters are unitary and normalized to be trivial on the diagonally embedded positive reals. For each $\pi \in \mathfrak{F}_m$, the associated $L$-function is defined by the Dirichlet series
\[
L(s,\pi) = \sum_{n=1}^{\infty}\frac{\lambda_\pi(n)}{n^s}.
\]
This series converges absolutely for $\Re(s) > 1$. We let $\mathbbm{1}$ denote the trivial representation, for which the associated $L$-function reduces to the Riemann zeta function, $L(s, \mathbbm{1}) = \zeta(s)$. We may therefore consider the function
\begin{align*}
    \frac{1}{L(s,\pi)} = \sum_{n=1}^\infty \frac{\mu_\pi(n)}{n^s},
\end{align*}
where the series also converges absolutely for $\Re(s) > 1$.
To analyze the behavior of $L(s, \pi)$ near the critical line $\Re(s) = 1/2$, we introduce a mollifier $M_{N, \pi}(s)$. For an integer $N \ge 2$, we consider the truncated Dirichlet polynomial
\begin{align*}
    M_{N,\pi}(s) = \sum_{n\le N}\frac{\mu_\pi(n)}{n^s}\left(1-\frac{\log n}{\log N}\right). 
\end{align*}
We then consider the mollified second moment over the interval $[T_1, T_2]$ as
\begin{align}
    I_{N,\pi}(T_1,T_2) = \int_{T_1}^{T_2}|M_{N,\pi}(\tfrac{1}{2}+it)L(\tfrac{1}{2}+it,\pi)|^2\, dt.\label{def:I(N,pi)}
\end{align}
The definition of $M_{N,\pi}(s)$ can be extended  to $M_{x,\pi}(s) $ for any $x>0$ by writing
\begin{align}\label{eq:Mlogx}
    M_{x,\pi}(s) \log x = \sum_{n\le x} \frac{\mu_\pi(n)}{n^s}\log (x/n),
\end{align}
with $M_{1,\pi}(s)=0$. Observe that $M_{x,\pi}(s) = 0$ for $0<x\le 1$, so the definition of $I_{N,\pi}(T_1,T_2)$ can be naturally extended to $I_{x,\pi}(T_1,T_2)$.

 The study of mollified moments has long been a central tool in understanding lower and upper bounds for the size of $L$-functions as well as the distribution of nontrivial zeros; see, for instance, \cite{MR3181494,MR4689334,MR2887604,MR3597392,MR4045193,radziwill, MR3404560, MR4002302}. In this setting, one normally takes $N=T^\theta$ for $\theta<1$. One critical aspect of obtaining good results is to ensure that $\theta$ is as large as possible. Levinson \cite{Levinson} found an asymptotic formula for the mollified second moment of $\zeta(s)$:
\[
\lim_{T\to \infty} \frac{I_{T^\theta, \mathbbm{1}}(0, T)}{T} = 1+\frac{1}{\theta}, \qquad 0 < \theta < 1/2,
\]
which he used to deduce that the proportion of zeros on the critical line, $\kappa$, satisfies $\kappa > 1/3$. Conrey \cite{Conrey} later extended this valid range to $\theta < 4/7$, proving that $\kappa > 2/5$ (see also \cite{MR4045193}).

While it was initially believed that such asymptotic formulas would not hold for $\theta > 1$, Farmer \cite{Farmer} produced a heuristic argument suggesting that the mean square remains bounded for every $\theta > 0$, an assertion now known as the ``$\theta = \infty$ conjecture''. Farmer showed that this conjecture implies $\kappa = 1$, meaning that $100\%$ of the nontrivial zeros of $\zeta(s)$ lie on the critical line. 

Building on these heuristics, Bettin and Gonek \cite{BettinGonek2017} provided a rigorous analytic proof that the $\theta = \infty$ conjecture is actually a sufficient condition for the Riemann Hypothesis. Specifically, they showed that for $\zeta(s)$, an upper bound of the form $I_{N, \mathbbm{1}}(0, T) \ll T^{1+\varepsilon}$ for all $N \le T^\theta$ with $\theta > 1$ is sufficient to force a zero-free region of the form $\Re(s) > \frac{1}{2}+\frac{1}{2\theta}$. If this bound holds for arbitrarily large $\theta$, the Riemann Hypothesis follows.

In this paper, we generalize the results of Bettin and Gonek \cite{BettinGonek2017} to the setting of general automorphic $L$-functions. 
The results of the present work show that the connection between arithmetic mollification and zero distribution is a universal phenomenon for automorphic $L$-functions. In addition, we develop a framework for deriving zero-free regions across broad families of automorphic representations. 

\begin{thm}\label{thm: main theorem-local}
Let $m \in \mathbb{N}$ and  $\pi\in\mathfrak{F}_m$. Let $\theta >0$, $\sigma\ge 1/2$, and $1/2 > \varepsilon> 0$ be fixed. There exists a constant $\Cl[abcon]{upper-condition} := \Cr{upper-condition}(\pi, \theta, \sigma, \varepsilon)  > 0$ with the following property: 

For any $T_2 \ge 1$, $0 \le T_1 \le T_2/2$, and $I_{y,\pi}(T_1, T_2)$  defined as in \eqref{def:I(N,pi)} 
for which
    \[
    \int_1^{T_2^{\theta}} I_{y,\pi}(T_1, T_2) \, dy \le \Cr{upper-condition} T_2^{2\sigma\theta-5-2m},
    \]
   the $L$-function $L(s,\pi)$  does not vanish in the rectangle defined by
    \[
    \Re(s) > \sigma+\varepsilon \quad \text{and} \quad T_1\le \Im(s) \le T_2.
    \]
\end{thm}

\begin{rek}
     In particular, for any $T \ge 1$, Theorem \ref{thm: main theorem-local} applies in the cases $(T_1,T_2)=(0,T)$ and $(T,2T)$. We also note that the assumptions in Theorem \ref{thm: main theorem-local} differ from those in \cite{BettinGonek2017}, where the assumption is an upper bound on $I_{y,\pi}(T_1, T_2)$ itself. This change represents a relaxation of the required conditions, as a pointwise bound is slightly stronger than necessary for the underlying arguments.
     This observation applies to the rest of our results as well.
\end{rek}

By strengthening the assumption on the mollified second moment, the same framework yields the non-vanishing of $L(s, \pi)$ in the entire half-plane $\Re(s) > \sigma$.
 \begin{thm}\label{thm: main theorem-all-T}
Let $m \in \mathbb{N}$ and $\pi\in\mathfrak{F}_m$. Let $\theta>0$ and $\sigma\ge 1/2$ be fixed. Suppose that 
    \begin{equation}\label{eq:condition-all T}
        \sup_{T \in [1,\infty)}  T^{-2\sigma\theta}\int_1^{T^{\theta}} I_{y,\pi}(0, T)\, dy <\infty.
    \end{equation}
Then  $L(s,\pi)$ does not vanish  in the half-plane $\Re(s) > \sigma$.
 \end{thm}



\begin{rek}\label{remark-Siegel-zero}

    It is worth mentioning that if $\sigma$ is sufficiently small, i.e. if there exists a fixed constant $\delta > 0$ such that the supremum condition holds for some $\sigma \le 1-\delta$, then the supremum condition in Theorem \ref{thm: main theorem-all-T} prohibits the existence of possible Landau--Siegel zeros. Following the proof of the theorem, one sees that the supremum condition forces every zero $\beta+i\gamma$ of $L(s,\pi)$ to satisfy $\beta\leq \sigma$. In fact, the existence of Landau--Siegel zeros would blow up the mollified moments and violate the supremum condition.
\end{rek}
As an immediate consequence of Theorem \ref{thm: main theorem-all-T}, we obtain Corollary \ref{cor:RH} by setting $\sigma = \tfrac{1}{2} + \varepsilon$ and allowing $\varepsilon \to 0^+$. This corollary serves as an analogue to the main result of Bettin and Gonek \cite{BettinGonek2017} when $\theta$ is arbitrarily large.
\begin{cor}\label{cor:RH}
Let $m \in \mathbb{N}$ and $\pi\in\mathfrak{F}_m$. Suppose there exists $\theta>0$ such that for every $\varepsilon>0$
  \[
  \sup_{T \in [1,\infty)}T^{-(1+\varepsilon)\theta}\int_{1}^{T^\theta}I_{y,\pi}(0, T) \, dy< \infty.
 \]
 Then the Riemann Hypothesis holds for $L(s,\pi)$. 
\end{cor}

\begin{rek}
   In the case of the Riemann zeta function $(\pi \cong \mathbbm{1})$, we note that the condition in Corollary \ref{cor:RH} is a quantitative analogue of the $\theta = \infty$ conjecture considered by Bettin and Gonek \cite{BettinGonek2017}. This supremum condition is expected to be more accessible in the regime where $\theta$ is large. Specifically, suppose, as in \cite[Theorem~1]{BettinGonek2017}, that $I_{N, \mathbbm{1}}(0,T) \ll_{\varepsilon'} T^{1+\varepsilon'}$ for all $2 \le N \le T^{\theta}$. Then for any $\varepsilon > 0$, we have
   \[
   T^{-(1+\varepsilon)\theta} \int_{1}^{T^\theta} I_{y,\mathbbm{1}}(0, T) \, dy \ll_{\varepsilon'} T^{1 + \theta(1+\varepsilon') - (1+\varepsilon)\theta} = T^{1 + \theta(\varepsilon' - \varepsilon)}.
\]
This expression is bounded uniformly for $T \ge 1$ provided that $\theta \ge (1 + \theta\varepsilon') / \varepsilon$. Consequently, for the supremum condition to be satisfied for arbitrarily small $\varepsilon > 0$, one must allow $\theta$ to be correspondingly large, consistent with the asymptotic nature of the original conjecture in \cite{BettinGonek2017}.
\end{rek}

\begin{rek}\label{remark-all-T}
    Alternatively, the same zero-free region in Theorem \ref{thm: main theorem-all-T} is implied by the condition
    \[
    \sup_{T \in [1,\infty)} \, T^{-2\sigma\theta + 3} \int_1^{T^{\theta}} I_{y,\pi}(T, 2T) \, dy < \infty.
    \]
    Consequently, the Riemann Hypothesis for $L(s,\pi)$ is implied by the existence of $\theta > 0$ such that for any $\varepsilon >0$,
    \[
    \sup_{T \in [1,\infty)} \, T^{-(1+\varepsilon)\theta+3} \int_1^{T^{\theta}} I_{y,\pi}(T, 2T) \, dy < \infty.
    \]
    The comparison between the conditions in Theorem \ref{thm: main theorem-all-T} and Corollary \ref{cor:RH} and the conditions in this remark suggests that shifting the analysis to a dyadic interval necessitates a stronger condition on the mollified moment. This phenomenon was also noted in \cite{BettinGonek2017}.
 \end{rek}

Building on the estimates established above, we extend our framework to provide the following results for families of cuspidal automorphic $L$-functions. Let $\mathcal{F}$ be a subset of $\mathfrak{F}_m$, not necessarily finite, and let $N(\sigma, T_1, T_2; \mathcal{F})$ denote the number of $L$-functions $L(s, \pi) \in \mathcal{F}$ possessing at least one zero $\rho_{\pi} = \beta_{\pi} + i\gamma_{\pi}$ satisfying $\beta_{\pi} \ge \sigma$ and $\gamma_{\pi} \in [T_1, T_2]$.


\begin{thm}\label{thm: main theorem-average}
Let $m \in \mathbb{N}$,  $\pi\in\mathfrak{F}_m$, and $\mathcal{F} \subset \mathfrak{F}_m$. Let $\theta > 0$ and $\sigma \ge 1/2$ be fixed. Suppose that
there exists some $\delta > 0$ such that uniformly for $T \ge 1$,
\begin{equation}\label{eq:condition-family}
    \sum_{\pi \in \mathcal{F}}\int_1^{T^{\theta}} I_{y,\pi}(0,T) \, dy \ll_{\delta} |\mathcal{F}|T^{2\sigma\theta-\delta}.
\end{equation}
Then for any $\varepsilon>0$,
\[
\frac{N\left(\sigma, 0, T; \mathcal{F}\right)}{|\mathcal{F}|} \ll_{\varepsilon} T^{-\delta+\varepsilon}.
\]
\end{thm}

\begin{rek}
   Taking $T \to \infty$ in the conclusion of Theorem \ref{thm: main theorem-average}, we obtain
   \[
   \lim_{T \to \infty} \frac{N\left(\sigma, 0, T; \mathcal{F}\right)}{|\mathcal{F}|} = 0.
\]
This suggests that for $100\%$ of the family $\mathcal{F}$, the associated $L$-functions do not vanish in the half-strip $\Re(s) > \sigma$, establishing a quasi-Riemann Hypothesis for the family. In the case where $\mathcal{F}$ is a fixed finite set, this asymptotic behavior implies that all $L$-functions in the family are simultaneously non-vanishing in the half-strip $\Re(s) > \sigma$.
\end{rek}

\section{Properties of automorphic \texorpdfstring{$L$-functions}{L-functions}}
In this section and throughout the remainder of the paper, we adopt the standard notation $f(z) = O_{\nu}(g(z))$, or equivalently $f(z) \ll_{\nu} g(z)$, to denote that $|f(z)| \le c|g(z)|$ for some effectively computable constant $c$ depending only on the parameter $\nu$. If no subscript is present, the implied constant is understood to be absolute. Furthermore, we let $\varepsilon > 0$ denote an arbitrarily small constant that may vary from line to line. 

   Let $m \ge 1$ be an integer,  $\mathbb{A}_{\mathbb{Q}}$ the ring of ad\`{e}les over $\mathbb{Q}$, and $\mathfrak{F}_{m}$ the set of all cuspidal automorphic representations of $\mathrm{GL}_{m}(\mathbb{A}_{\mathbb{Q}})$ whose central characters are unitary and normalized to be trivial on the diagonally embedded positive reals. For any $\pi \in\mathfrak{F}_{m}$, there exists a smooth admissible representation $\pi_{v}$ of $\GL_{m}(\mathbb{Q}_v)$ at each place $v$ of $\mathbb{Q}$ such that $\pi$ decomposes as the restricted tensor product $\otimes_v \pi_{v}$. 

Given $\pi \in\mathfrak{F}_{m}$, we denote by  $\mathfrak{q}_\pi$ its arithmetic conductor.  At a non-archimedean place $v$ corresponding to a prime $p$, the local $L$-function $L(s,\pi_p)$ is defined through the Satake parameters $\alpha_{1,\pi}(p),\alpha_{2,\pi}(p),\cdots,\alpha_{m,\pi}(p)$ as 
\[
L(s,\pi_p)=\prod_{r=1}^m(1-\alpha_{r,\pi}(p)p^{-s})^{-1}=\sum_{k=0}^\infty \frac{\lambda_\pi(p^k)}{p^{ks}},
\]
where the coefficient $\lambda_\pi(p^k)$ is given by
\[
\lambda_\pi(p^k) = \sum_{\substack{\ell_1+\ell_2+\cdots+\ell_m=k\\\ell_i\in\N\cup\{0\}}}\alpha_{1,\pi}(p)^{\ell_1}\alpha_{2,\pi}(p)^{\ell_2}\cdots\alpha_{m,\pi}(p)^{\ell_m}. 
\]
If $p\nmid \mathfrak{q}_\pi$, then for all $1\le r\le m$, we have $\alpha_{r,\pi}(p)\neq 0$. If $p\mid \mathfrak{q}_\pi$, then there might exist $r$ such that $\alpha_{r,\pi}(p)= 0$.
The standard $L$-function $L(s,\pi)$ associated to $\pi$ is defined by the Euler product and Dirichlet series
    \begin{equation}\label{def:dirichlet series for Lu(s)}
        L(s,\pi)=\prod_{p} L(s,\pi_{p})=\sum_{n=1}^{\infty}\frac{\lambda_{\pi}(n)}{n^s},
    \end{equation}
both of which converge absolutely when $\Re(s)>1$. 

 At the archimedean place of $\mathbb{Q}$, the local $L$-function $ L(s,\pi_{\infty})$ is defined using $m$ spectral parameters $\nu_{\pi}(r)\in\mathbb{C}$:
    \[
        L(s,\pi_{\infty}) = \pi^{-\frac{ms}{2}}\prod_{r=1}^{m}\Gamma\Big(\frac{s+\nu_{\pi}(r)}{2}\Big).
    \]
    Let $\widetilde{\pi}\in\mathfrak{F}_{m}$ be the contragredient representation of $\pi$. We have $\mathfrak{q}_{\pi}=\mathfrak{q}_{\widetilde{\pi}}$, and 
    \[
        \{\alpha_{r,\widetilde{\pi}}(p)\}=\{\overline{ \alpha_{r,\pi}(p)}\},\qquad \{\nu_{\widetilde{\pi}}(r)\}=\{\overline{\nu_{\pi}(r)}\}.
    \]
    Let $\omega_{\pi}$ be the order of the pole of $L(s,\pi)$ at $s=1$: this is $0$ unless $m = 1$ and $\pi =\mathbbm{1}$.
    The completed $L$-function
    \[
        \Lambda(s,\pi) = (s(1-s))^{\omega_{\pi}}\mathfrak{q}_{\pi}^{s/2}L(s,\pi)L(s,\pi_{\infty})
    \]
    is entire of order $1$, and there exists a complex number $W(\pi)$ of modulus $1$ such that for all $s\in\mathbb{C}$, the functional equation $\Lambda(s,\pi)=W(\pi)\Lambda(1-s,\widetilde{\pi})$ holds. The analytic conductor of $\pi$ \cite{IS} is given by
    \begin{equation}
	\label{eqn:analytic_conductor_def for Lu(s)}
           \mathfrak{C}(s, \pi)\coloneqq \mathfrak{q}_{\pi}\prod_{r=1}^{m} (3+|s+\nu_{\pi}(r)|),\qquad \mathfrak{C}(\pi)\coloneqq \mathfrak{C}(0,\pi).
    \end{equation}

    The Riemann Hypothesis (RH) states that all zeros of the completed $L$ function $\Lambda(s,\pi)$ are on the critical line $\Re(s) = \tfrac{1}{2}$. The Generalized Ramanujan--Petersson Conjecture (GRC) states that the Satake parameters $\alpha_{r,\pi}(p)$  satisfy $|\alpha_{r,\pi}(p)|=1$ for all but a finite number of primes $p$. In general, this conjecture is open. Towards GRC, Luo, Rudnick, and Sarnak \cite{LRS} and M\"uller and Speh \cite{MS} established the existence of a constant $\delta_{m}\in[0,\frac{1}{2}-\frac{1}{m^2+1}]$ such that we have
    \begin{equation}
	\label{eqn:LRS_finite}
            |\alpha_{r,\pi}(p)|\le  p^{\delta_{m}}\qquad\text{ and }\qquad\Re(\nu_{\pi}(r))\ge -\delta_{m},
    \end{equation}
    and GRC predicts that one may take $\delta_{m}=0$ in \eqref{eqn:LRS_finite}.


\section{Proofs of the main theorems}
The proofs for Theorems \ref{thm: main theorem-local}, \ref{thm: main theorem-all-T}, and \ref{thm: main theorem-average} rely on a similar framework, and thus we will only present the proof of Theorem \ref{thm: main theorem-local} in detail, and point out the differences in the proofs for Theorem \ref{thm: main theorem-all-T} and Theorem \ref{thm: main theorem-average}.

\subsection{Initial setup for the proofs}
For $c>0$, we recall the identity
\begin{align}
\frac{1}{2\pi i}\int_{c-i\infty}^{c+i\infty} y^z \ \frac{dz}{z^2}
=
\begin{cases}
\log y, & y>1,\\[6pt]
0, & 0<y<1.
\end{cases}\label{eq:log identity}
\end{align}
Since we can write the right side of \eqref{eq:Mlogx} as 
\[
 \sum_{n \le x} \frac{\mu_{\pi}(n)}{n^s}\log(x/n) = \sum_{n=1}^{\infty}\frac{\mu_{\pi}(n)}{n^s} \frac{1}{2\pi i}\int_{c-i\infty}^{c+i\infty} (x/n)^z \ \frac{dz}{z^2}=\frac{1}{2\pi i} \int_{1-i\infty}^{1+i\infty}x^z\sum_{n=1}^{\infty}\frac{\mu_\pi(n)}{n^{\tfrac{1}{2}+it+z}}\ \frac{dz}{z^2},
\]
we obtain
\begin{equation}
     M_{x,\pi}(\tfrac{1}{2}+it) \log x 
     = \frac{1}{2\pi i} \int_{1-i\infty}^{1+i\infty}\frac{x^z}{L(\tfrac{1}{2}+it+z,\pi)} \ \frac{dz}{z^2}.\label{eq:M(1/2+it)log x}
\end{equation}

We define $H_t(w,\pi)$ as the Mellin transform of $M_{x,\pi}(\tfrac{1}{2}+it) \log x$:
\begin{align*}
    H_t(w,\pi) := \int_1^\infty  M_{x,\pi}(\tfrac{1}{2}+it) (\log x) x^{-w}\, dx. 
\end{align*}
For any $\varepsilon >0$, we observe that
\begin{align*}
   |M_{x,\pi}(\tfrac{1}{2}+it) |\le \sum_{n\le x}\frac{|\mu_\pi(n)|}{n^{\frac{1}{2}}}\ll_{\varepsilon} \sum_{n\le x}n^{-\frac{1}{2}+\delta_m+\varepsilon},
\end{align*}
using \eqref{eqn:LRS_finite} and the fact that $|\mu_{\pi}(n)| \le \tau_m(n)n^{\delta_m}$. Here, $\tau_m(n)$ denotes the generalized divisor function, which counts the number of ways to represent $n$ as an ordered product of $m$ positive integers. It follows that
\begin{align*}
   |M_{x,\pi}(\tfrac{1}{2}+it) (\log x) x^{-w}| \ll_\varepsilon x^{\frac{1}{2}+\delta_m+\varepsilon-\Re(w)}.
\end{align*}
Since $\delta_m \le \frac{1}{2}-\frac{1}{m^2+1}$, the integral $H_t(w,\pi)$ converges absolutely for $\Re(w)\ge 2$. Substituting \eqref{eq:M(1/2+it)log x} into the definition of $H_t(w,\pi)$ and computing the residue at $z=w-1$, we obtain for $\Re(w)\ge 2$,
\begin{align*}
    H_t(w,\pi)  = \frac{1}{(w-1)^2L(w-\tfrac{1}{2}+it,\pi)}.
\end{align*}

\subsection{Proof of Theorem \ref{thm: main theorem-local}}

We recall the assumption in Theorem \ref{thm: main theorem-local} and proceed by contradiction. Consider a fixed pair $T_2 \ge 1$ and $0 \le T_1 \le T_2/2$. Suppose for the sake of contradiction that there exists $0 < \varepsilon_0 < 1/2$  such that for every $\Cr{upper-condition} > 0$, the hypothesis
\[
\int_{1}^{T_2^\theta} I_{y,\pi}(T_1, T_2) \, dy \le \Cr{upper-condition} T_2^{2\sigma\theta-5-2m}
\]
holds, yet there exists a zero $\rho_0 = \beta_0 + i\gamma_0$ of $L(s, \pi)$ such that $\beta_0 > \sigma+\varepsilon_0$ and $ \gamma_0 \in [T_1, T_2]$. We shall show that such an assumption leads to a contradiction if  $\Cr{upper-condition}$ is  sufficiently small.

Define
\[
G_t(w,\pi):= \frac{(w-1)^2(w-\tfrac{3}{2}+it)^{\delta_{\pi = \mathbbm{1}}}L(w-\tfrac{1}{2}+it,\pi)}{(w+1)^2(w-\tfrac{1}{2}+it-\rho_0)(w+it+1)^{1+m+\delta_{\pi = \mathbbm{1}}}},
\]
where $\delta_{\pi = \mathbbm{1}}=1$ is the Kronecker delta function for $L(s,\pi)=\zeta(s)$.
The function $ G_t(w,\pi)$ is holomorphic for $\Re(w)\ge 0$. When $\Re(w) > \tfrac{3}{2}$, we obtain
\[
\frac{(w-\tfrac{3}{2}+it)^{\delta_{\pi = \mathbbm{1}}}L(w-\tfrac{1}{2}+it, \pi)}{(|w+it|+1)^{\delta_{\pi = \mathbbm{1}}}} \ll 1.
\]
When $\frac{1}{2} \le \Re(w) \le \frac{3}{2}$, the Phragmen-Lindel\"of principle implies that
\begin{align*}
    \frac{(w-\tfrac{3}{2}+it)^{\delta_{\pi = \mathbbm{1}}}L(w-\tfrac{1}{2}+it,\pi)}{(|w+it|+1)^{\delta_{\pi = \mathbbm{1}}}} &\ll_{\varepsilon} \mathfrak{C}(w-\tfrac{1}{2}+it,\pi)^{\frac{3/2-\Re(w)}{2}+\varepsilon}\\
    &\ll_{\pi, \varepsilon} (|w+it|+1)^{\frac{m}{2}+\varepsilon}.
\end{align*}
Lastly, when $0 \le \Re(w) \le \frac{1}{2}$, we approximate $L(-w-\frac{1}{2}+it, \pi)$ using the functional equation and Stirling's formula $|\Gamma(\sigma + it)| \ll  |t|^{\sigma - 1/2} e^{-\pi |t| / 2}$: 
\begin{equation*}
    \begin{aligned}
        L(w-\tfrac{1}{2}+it, \pi) &\ll_{\pi} \left|\frac{ \prod_{r=1}^m \Gamma\left(\frac{3/2 - w-it+ \nu_{\widetilde{\pi}}(r)}{2}\right)}{\prod_{r=1}^m \Gamma\left(\frac{-1/2 +w+it+ \nu_{\pi}(r)}{2}\right)} \right| |L(\tfrac{3}{2}-w-it, \widetilde{\pi})|\\
        &\ll_{\pi} \frac{(|w+it|+1)^{m(3/2-\Re(w))/2}}{(|w+it|+1)^{m(-1/2+\Re(w))/2}} \ll_{\pi}(|w+it|+1)^m. 
    \end{aligned}
\end{equation*}

Combining the preceding estimates, we find that $G_t(w)$ is holomorphic in the half-plane $\Re(w) \ge 0$ and satisfies the bound
\begin{align*}
G_t(w,\pi) 
&\ll \frac{(|w+it|+1)^{m+\varepsilon}}{(w-\tfrac{1}{2}+it-\rho_0)(|w+it|+1)^{{1+{m}}}}\\
&\ll_{\varepsilon} (|w+it|+1)^{-2+\varepsilon}.
\end{align*}
Setting 
\[
g_t(u) = \frac{1}{2\pi i}\int_{3-i\infty}^{3+i\infty} G_t(w,\pi)u^{-w}\, dw
\]
for $u > 0$, we have 
\begin{equation}\label{eq:g_t(u)}
g_t(u) = \begin{cases}
    0 \qquad &\text{ if } u > 1,\\
    O(1) \qquad &\text{ if } 0 \le u \le 1, 
\end{cases}
\end{equation}
as can be seen by moving the line of integration to $\Re(w) = \infty$ when $u > 1$, and to $\Re(w) = 0$ when $0 \le u \le 1$.

Now consider the integral, 
\begin{equation}\label{eq:J_t}
\begin{aligned}
    J_t(x) &:= \frac{1}{2\pi i} \int_{3-i\infty}^{3+i\infty}G_t(w,\pi)H_t(w,\pi)x^w \, dw\\
    &= \frac{1}{2\pi i} \int_{3-i\infty}^{3+\infty} \frac{x^w(w-\tfrac{3}{2}+it)^{\delta_{\pi = \mathbbm{1}}}}{(w+1)^2(w-\tfrac{1}{2}+it-\rho_0)(w+it+1)^{1+m+\delta_{\pi = \mathbbm{1}}}}\, dw,
\end{aligned}
\end{equation}
where from this point on  we assume that $x \ge 2$.

On the one hand, we recall that the Mellin transform of a multiplicative convolution is the product of the constituent Mellin transforms. Specifically, for functions $f$ and $g$, let
\[
(f \star g)(x) := \int_{0}^{\infty} f(y) g(x/y) \, \frac{dy}{y}.
\]
Then the transform satisfies
\[
\int_{0}^{\infty} (f \star g)(x) x^{s-1} dx = \left( \int_{0}^{\infty} f(x) x^{s-1} \, dx \right) \left( \int_{0}^{\infty} g(x) \ x^{s-1} dx \right).
\]
Since $M_{y,\pi}(\tfrac{1}{2}+it)\log y$ vanishes for $y \le 1$, the convolution simplifies to
\[
J_t(x) = \int_1^{\infty} M_{y,\pi}(\tfrac{1}{2}+it)(\log y) g_t(x/y)\ \frac{dy}{y}.
\]
This can be seen as follows: identifying $H_t(w,\pi)$ as the Mellin transform of $M_{y,\pi}(\tfrac{1}{2}+it)\log y$ and $G_t(w)$ as the transform of $g_t(u)$, it follows that $J_t(x)$ is the Mellin inversion of the product $G_t(w) H_t(w,\pi)$.
Therefore, by \eqref{eq:g_t(u)}, we have for $x \ge 2$,
\begin{equation}\label{eq:bound-Jt}
    J_t(x) \ll \int_{1}^{x} |M_{y,\pi}(\tfrac{1}{2}+it)| (\log y)\, dy.
\end{equation}

 On the other hand, moving the line of integration in \eqref{eq:J_t} to $\Re(w)=0$, we see that 
\begin{equation}\label{eq:J-t-push-contour}
J_t(x) = \frac{1}{2\pi i}\int_{-i\infty}^{i\infty} G_t(w,\pi)H_t(w,\pi)x^wdw + \frac{x^{\rho_0+\frac{1}{2}-it}(\rho_0-1)^{\delta_{\pi = \mathbbm{1}}}}{(\rho_0-it+\frac{3}{2})^2(\rho_0+\frac{3}{2})^{1+m+\delta_{\pi = \mathbbm{1}}}}
\end{equation}
by picking up the residue from the pole at $w=\rho_0+\tfrac{1}{2}-it$. The integral on the right is $O(1)$ since $G_t(w,\pi)H_t(w,\pi) \ll (1+|w|)^{-2}$ for $\Re(w)= 0$. Thus, from \eqref{eq:bound-Jt} and \eqref{eq:J-t-push-contour} we deduce  
\begin{align}
\frac{x^{\beta_0+\frac{1}{2}}}{(|\gamma_0-t|+1)^2(|\gamma_0|+1)^{1+m}}+1 \ll \int_{1}^{x} |M_{y,\pi}(\tfrac{1}{2}+it)| \log y\, dy.\label{eq:main change in proof 1}
\end{align}
Squaring both sides of \eqref{eq:main change in proof 1} and applying the Cauchy--Schwarz inequality, we obtain for $x \ge 2$
\[
\frac{x^{2\beta_0}}{(|\gamma_0-t|+1)^4(|\gamma_0|+1)^{2+2m}}+\frac{1}{x} \ll \int_{1}^{x} |M_{y,\pi}(1/2+it)|^2 (\log y)^2\, dy.
\]
We note that $\gamma_0 \in [T_1, T_2]$. Multiplying both sides by $|L(\frac{1}{2}+it,\pi)|^2$ and integrating with respect to $t$ over the interval $[T_1, T_2]$, we obtain 
\begin{equation}\label{eq:multiplying-L}
    \begin{aligned}
        &\int_{T_1}^{T_2}|L(\tfrac{1}{2}+it, \pi)|^2 \left(\frac{x^{2\beta_0}}{(|\gamma_0-t|+1)^4(|\gamma_0|+1)^{2+2m}}+\frac{1}{x}\right) dt \\
        &\ll \int_{1}^{x} (\log y)^2 \int_{T_1}^{T_2}|M_{y,\pi}(\tfrac{1}{2}+it)L(\tfrac{1}{2}+it,\pi)|^2 \, dt\, dy\\
        &\le (\log x)^2 \int_1^{x} I_{y,\pi}(T_1, T_2)\, dy.
    \end{aligned}
\end{equation}

By applying \cite[Theorem 1.1]{ji2009} dyadically, we have for any $T_2 \ge 1$ and $0 \le T_1 \le T_2/2$ 
\begin{equation}\label{eq:lowerbound}
    \int_{T_1}^{T_2} |L(\tfrac{1}{2}+it, \pi)|^2 dt \gg T_2.
\end{equation}
We now approximate the left side of \eqref{eq:multiplying-L}. We first note that
\begin{equation}\label{eq:1/x contributaion}
        \int_{T_1}^{T_2}|L(\tfrac{1}{2}+it, \pi)|^2 \frac{1}{x}\,dt \gg  \frac{T_2}{x}.
\end{equation}
Defining $J(T) = \int_0^{T}|L(\tfrac{1}{2}+it, \pi)|^2\,dt$, by partial summation, we have
\begin{equation}\label{eq:main change in proof 2}
    \begin{aligned}
        &\int_{T_1}^{T_2}|L(\tfrac{1}{2}+it, \pi)|^2 \frac{x^{2\beta_0}}{(|\gamma_0-t|+1)^4(|\gamma_0|+1)^{2+2m}}\,dt
        \gg  \frac{x^{2\beta_0}}{T_2^{5+2m}}.
    \end{aligned}
\end{equation}

Putting together \eqref{eq:1/x contributaion} and \eqref{eq:main change in proof 2} in \eqref{eq:multiplying-L}, we see that there exists an effectively computable constant $\Cl[abcon]{lower} := \Cr{lower}(\pi) > 0$ such that
\begin{equation}\label{eq:optimization}
   \Cr{lower}\left(\frac{x^{2\beta_0}}{T_2^{5+2m}} + \frac{T_2}{x}\right) \le (\log x)^2 \int_1^{x} I_{y,\pi}(T_1, T_2) dy. 
\end{equation}
Take $x = T_2^\theta$. By assumption,  for every $\Cr{upper-condition} > 0$ we have
\[
    \log^2(T_2^{\theta}) \int_1^{T_2^\theta} I_{y,\pi}(T_1, T_2) dy \le \Cr{upper-condition} \log^2 (T_2^{\theta})  T_2^{2\sigma\theta-5-2m}.
 \]
Substituting in \eqref{eq:optimization}, we get 
\begin{equation}\label{eq:comparison}
     \Cr{lower}\left(\frac{T_2^{2\beta_0\theta}}{T_2^{5+2m}} + \frac{T_2}{x}\right) \le \Cr{upper-condition} \log^2(T_2^{\theta}) T_2^{2\sigma\theta-5-2m}.
\end{equation}
Since $\varepsilon_0 > 0$, for all $T_2 \ge 1$, there exists an effectively computable constant $\Cl[abcon]{log}:= \Cr{log}(\theta, \varepsilon_0) >0$ such that 
\[\log^2 (T_2^\theta) \le  \Cr{log} T_2^{2\theta\varepsilon_0}. 
\]

We then obtain from \eqref{eq:comparison} 
\[
\frac{\Cr{lower}}{\Cr{upper-condition}\Cr{log}} \le T_2^{-2(\beta_0-\sigma-\varepsilon_0)\theta} \le 1,
\]
since $\beta_0 > \sigma+\varepsilon_0$, $\theta > 0$, and $T_2 \ge 1$.
Thus, by choosing $\Cr{upper-condition}$ to be sufficiently small, specifically such that $0 < \Cr{upper-condition} < \Cr{lower}/\Cr{log}$, the left side becomes strictly greater than $1$, yielding a contradiction.

\subsection{Proof of Theorem \ref{thm: main theorem-all-T}} 
Unlike the proof of Theorem \ref{thm: main theorem-local}, which relies on a proof by contradiction, the proof of Theorem \ref{thm: main theorem-all-T} proceeds via a direct asymptotic comparison. Suppose that the condition \eqref{eq:condition-all T} holds true and consider a nontrivial zero $\rho_0 = \beta_0 + i\gamma_0$ of $L(s, \pi)$. We shall show that $\beta_0 \le \sigma$, as claimed.

The proof of Theorem \ref{thm: main theorem-all-T} parallels the argument for Theorem \ref{thm: main theorem-local}, and as such, we highlight only the necessary modifications. In this setting, the inequality in \eqref{eq:multiplying-L} is now changed to 
\begin{equation}\label{eq:changed 1-all T}
\begin{aligned}
\int_{0}^{T}|L(\tfrac{1}{2}+it, \pi)|^2 \left(\frac{x^{2\beta_0}}{(|\gamma_0-t|+1)^4(|\gamma_0|+1)^{2+2m}}+\frac{1}{x}\right) dt
\le (\log x)^2 \int_1^{x} I_{y,\pi}(0, T)\, dy.
\end{aligned}
\end{equation}
Using \eqref{eq:lowerbound} and following the calculation in \eqref{eq:1/x contributaion} and \eqref{eq:main change in proof 2}, we see that   
\begin{equation}\label{eq:changed optimization}
 \int_{0}^{T}|L(\tfrac{1}{2}+it, \pi)|^2 \left(\frac{x^{2\beta_0}}{(|\gamma_0-t|+1)^4(|\gamma_0|+1)^{2+2m}}+\frac{1}{x}\right) dt  \gg x^{2\beta_0} + \frac{T}{x}. 
\end{equation}
Substituting the lower bound \eqref{eq:changed optimization} into \eqref{eq:changed 1-all T} and setting $x = T^{\theta}$, we obtain
 \[
   T^{2\beta_0\theta} + T^{1-\theta} \ll   \log^2 (T^{\theta})\int_{1}^{T^\theta}I_{y,\pi}(0, T) \, dy.
 \]
Therefore, from the assumption \eqref{eq:condition-all T}, we see that for all $T \ge 1$ and any $\varepsilon > 0$,
\[
 T^{2\beta_0\theta} \ll_{\varepsilon} T^{2\sigma\theta + \varepsilon}.
\]
Taking $T \to \infty$ and $\varepsilon \to 0^{+}$, we obtain $\beta_0 \le \sigma$ as desired.

The extra factor $T^3$ in Remark \ref{remark-all-T} accounts for the lower bound on the dyadic integral version of \eqref{eq:changed 1-all T}. Specifically, we have
\[
\int_{T}^{2T}|L(\tfrac{1}{2}+it, \pi)|^2 \left(\frac{x^{2\beta_0}}{(|\gamma_0-t|+1)^4(|\gamma_0|+1)^{2+2m}}+\frac{1}{x}\right) dt \gg \frac{x^{2\beta_0}}{T^3} + \frac{T}{x},
\]
where the $T^{-3}$ decay in the first term mandates the additional power of $T$ in the conditions of Remark \ref{remark-all-T}.

\subsection{Proof of Theorem \ref{thm: main theorem-average}}
Similarly to the proof of Theorem \ref{thm: main theorem-all-T}, the proof of Theorem \ref{thm: main theorem-average} proceeds via a direct asymptotic comparison. Suppose that the condition \eqref{eq:condition-family} holds. We first define $\mathcal{B}(Q)$ as the set of $L$-functions in $\mathcal{F}$ that have at least one zero $\rho_\pi = \beta_\pi + i\gamma_\pi$ in the rectangle of consideration, namely, 
\[
\mathcal{B} = \{ \pi \in \mathcal{F} : \text{ there exists } \rho_\pi \text{ with } \beta_\pi \ge \sigma \text{ and } \gamma_\pi \in [0, T] \}.
\]
Similarly to \eqref{eq:changed 1-all T}, for each $\pi \in \mathcal{B}$, we have
\begin{equation}\label{eq:change-ave-1}
    \int_{0}^{T}|L(\tfrac{1}{2}+it, \pi)|^2 \left(\frac{x^{2\beta_{\pi}}}{(|\gamma_{\pi}-t|+1)^4(|\gamma_{\pi}|+1)^{2+2m}}+\frac{1}{x}\right) dt
\le (\log x)^2 \int_1^{x} I_{y,\pi}(0, T)\, dy.
\end{equation}
Using \eqref{eq:lowerbound} and following the calculation in \eqref{eq:1/x contributaion} and \eqref{eq:main change in proof 2}, we see that 
\begin{equation}\label{eq:lower-change-ave}
   \int_{0}^{T}|L(\tfrac{1}{2}+it, \pi)|^2 \left(\frac{x^{\beta_{\pi}}}{(|\gamma_{\pi}-t|+1)^4(|\gamma_{\pi}|+1)^{2+2m}}+\frac{1}{x}\right) dt  \gg x^{2\beta_{\pi}} + \frac{T}{x}.  
\end{equation}
Substituting the lower bound \eqref{eq:lower-change-ave} into \eqref{eq:change-ave-1}, setting $x = T^{\theta}$, and summing both sides over the set $\mathcal{B}$, we obtain
\begin{equation}\label{eq:lowerbound-ave}
       \sum_{\pi \in \mathcal{B}}(T^{2\beta_{\pi}\theta} + T^{1-\theta}) \ll   \sum_{\pi \in \mathcal{B}}\log^2 (T^{\theta})\int_{1}^{T^\theta}I_{y,\pi}(0, T) \, dy.
\end{equation}
From the assumption \eqref{eq:condition-family} and the fact that the sum over the bad members cannot exceed the sum over the entire family, there exists $\delta >0$ such that for any $\varepsilon> 0$,
\begin{equation}\label{eq:upperbound-ave}
    \sum_{\pi \in \mathcal{B}}\log^2 (T^{\theta})\int_{1}^{T^\theta}I_{y,\pi}(0, T) \, dy \ll_{\varepsilon} |\mathcal{F}|T^{2\sigma\theta-\delta+\varepsilon}
\end{equation}
Putting together \eqref{eq:lowerbound-ave} and \eqref{eq:upperbound-ave} and using the fact that $\beta_{\pi} \ge \sigma$, we obtain
\[
|\mathcal{B}|T^{2\sigma\theta} \ll_{\varepsilon} |\mathcal{F}|T^{2\sigma\theta-\delta+\varepsilon}.
\]
Since $|\mathcal{B}| = N(\sigma,0,T;\mathcal{F})$, the conclusion of Theorem \ref{thm: main theorem-average} follows immediately.

\section*{Funding}
A.D. and N.W. are supported by the Shaff--Andrews Fellowship, Department of Mathematics, University of Illinois Urbana-Champaign.
\bibliographystyle{abbrv}
\bibliography{theta_infty}
\end{document}